\documentclass[11pt,reqno,fleqn]{article}

\usepackage{amsmath,amsthm,amssymb}
\usepackage{mathpazo}
\usepackage{graphicx}
\usepackage{hyperref}
\usepackage{amsrefs}

\newtheorem{theorem}{Theorem}
\newtheorem{prop}[theorem]{Proposition}
\newtheorem{lemma}[theorem]{Lemma}
\newtheorem{conjecture}[theorem]{Conjecture}
{
\theoremstyle{remark}
\newtheorem*{definition}{Definition}
\newtheorem*{example}{Example}
\newtheorem*{remark}{Remark}
}

\newcommand{\MRnumber}[1]{~\href{http://www.ams.org/mathscinet-getitem?mr=#1}{{\bf MR~#1}}}

\newcommand{\RR}{{\mathbb R}}
\newcommand{\Z}{{\mathbb Z}}

\newcommand{\ZZ}{\ensuremath{\mathbb Z}}

\newcommand{\Prob}[1]{\ensuremath{{\mathbb P}\left[ #1 \right]}}

\renewcommand{\mod}[1]{{\ifmmode\text{\rm\ (mod~$#1$)}\else\discretionary{}{}{\hbox{ }}\rm(mod~$#1$)\fi}}

\newcommand{\sd}{sum-dominant}
\newcommand{\dd}{difference-dominant}
\newcommand{\sdb}{sum-difference-balanced}

\newcommand{\nset}{\{0, 1, 2, \dots, n-1\}}
\newcommand{\smallnset}{\{0,\dots,n-1\}}

\title{Many Sets Have More Sums Than Differences}
\author{Greg Martin \\ {\tt gerg@math.ubc.ca} \\ University of British Columbia \\ and \\
    Kevin O'Bryant \\ {\tt kevin@member.ams.org} \\ City University of New York, Staten Island}

\begin{document}
\maketitle \pagestyle{myheadings}\thispagestyle{empty} \markboth{GREG MARTIN AND KEVIN
O'BRYANT}{MANY SETS HAVE MORE SUMS THAN DIFFERENCES}

\section{Introduction}

As addition is commutative but subtraction is not, the set of sums
    \[
    S+S :=\{s_1+s_2 \colon s_i\in S\}
    \]
of a finite set $S$ is predisposed to be smaller than the set of differences
    \[
    S-S :=\{s_1-s_2 \colon s_i \in S\}.
    \]
As Nathanson~\cite{arxiv:math.NT/0604340} wrote:
\begin{quote}
``Even though there exist sets $A$ that have more sums than
differences, such sets should be rare, and it must be true with the right way of counting that the
vast majority of sets satisfies $|A-A|>|A+A|$.''
\end{quote}
Following this reasoning, one would suspect that a vanishingly small proportion of the $2^n$ subsets of
$\nset$ have more sums than differences. Our purpose, however, is to show that this is not the case. The following terminology will be used throughout this article:

\begin{definition}
A finite set $S$ is {\em\dd} if $|S-S|>|S+S|$, {\em\sd} if $|S+S|>|S-S|$, and {\em\sdb} if
$|S+S|=|S-S|$.
\end{definition}

\noindent
Nathanson~\cite{arxiv:math.NT/0608148} calls \sd\ sets ``MSTD'' sets, short
for ``More Sums Than Differences''. We refer the reader to~\cites{arxiv:math.NT/0604340,arxiv:math.NT/0608148} for the history of this problem.

Our main theorem shows that, perhaps contrary to intuition, all three types of set in the above definition are ubiquitous.

\begin{theorem}
Let $P$ be any arithmetic progression of length $n$. A positive proportion of the subsets of $P$
are \dd, a positive proportion are \sd, and a positive proportion are \sdb. More precisely, there
exists $c>0$ such that for all $n\ge15$,
\begin{align*}
\#\big\{ S\subseteq P &: S\text{ is \dd} \big\} > c2^n, \\
\#\big\{ S\subseteq P &: S\text{ is \sd} \big\} > c2^n, \\
\#\big\{ S\subseteq P &: S\text{ is \sdb} \big\} > c2^n.
\end{align*}
\label{positive.proportion.thm}
\end{theorem}

\noindent We observe that the sizes of $S+S$ and $S-S$ are invariant under translation and dilation of
$S$, so that without loss of generality we can restrict our attention to \mbox{$P=\nset$}.

The following examples show that none of the three categories is empty for $n\ge15$:
\begin{example}
The set $S=\{0,1,3\}$ has $S+S=\{0,1,2,3,4,6\}$ and $S-S=\{-3,-2,-1,0,1,2,3\}$; therefore $S$
is \dd, since $|S-S|=7>6=|S+S|$.
\end{example}
\begin{example}
The set $S=\{0,2,3,4,7,11,12,14\}$ has $S+S=\{0,\dots,28\} \setminus \{1, 20, 27\}$ and
$S-S=\{-14,\dots,14\} \setminus \{-13, -6, 6, 13\}$; therefore $S$ is \sd, since
$|S+S|=26>25=|S-S|$.
\end{example}
\begin{example}
A set $S$ is symmetric if $S=a^\ast -S$ for some $a^\ast\in\RR$. Any symmetric set has
$S+S=S+(a^\ast-S)=a^\ast+(S-S)$; therefore symmetric sets are \sdb. In particular, any interval or
arithmetic progression is \sdb.
\end{example}

The idea behind Theorem~\ref{positive.proportion.thm} is the following. Most subsets of $\nset$
have about $n/2$ elements; call our typical subset $S$. Each $k\in\{0,1,2,\dots,2n-2\}$ has, on
average, roughly $n/4-|n-k|/4$ representations as a sum of two elements of $S$. Not only is this positive,
it is quite large except when $k$ is near 0 or $2n-2$. Similarly, each nonzero
$k\in \{-(n-1),\dots,n-1\}$ has, on average, roughly $n/4-|k|/4$ representations as a difference of two
elements of $S$. Not only is this positive, it is quite large except when $|k|$ is near $n-1$.
Putting these together, the sizes of the sumset and difference set are predominantly affected by
the elements of $S$ that are near 0 or near $n$. If we choose the ``fringe'' of $S$ cleverly, the
middle of $S$ will become largely irrelevant.

This philosophy suggests the following conjecture; see Section~\ref{section.conjectures} for a more
refined conjecture.

\begin{conjecture}
Let $P$ be any arithmetic progression with length $n$. The limiting proportions
\begin{align*}
\rho_- &= \lim_{n\to\infty} 2^{-n} \,\#\big\{ S\subseteq P : S\text{ is \dd} \big\} \\
\rho_+ &= \lim_{n\to\infty} 2^{-n} \,\#\big\{ S\subseteq P : S\text{ is \sd} \big\} \\
\rho_= &= \lim_{n\to\infty} 2^{-n} \,\#\big\{ S\subseteq P : S\text{ is \sdb} \big\}
\end{align*}
all exist and are positive.
\label{main.conjecture}
\end{conjecture}

The following result, on the other hand, supports Nathanson's instinct as quoted above, with one interpretation of ``the right way'' and a suitably humble understanding of ``vast''. Theorem~\ref{theorem.average.value} is proved in Section~\ref{section.average.values}.
\begin{theorem}
Let $P$ be any arithmetic progression with length $n$. On average, the difference set of a subset
of $P$ has 4 more elements than its sumset. More precisely,
    \begin{align*}
    \frac1{2^n} \sum_{S\subseteq P} |S-S| &\sim 2n-7, \\
    \frac1{2^n} \sum_{S\subseteq P} |S+S| &\sim 2n-11.
    \end{align*}
\label{theorem.average.value}
\end{theorem}

Nathanson~\cite{arxiv:math.NT/0608148} asks for the possible values of $|A+A|-|A-A|$. We show by
construction in Section \ref{section.construction} that the range of $|A+A|-|A-A|$ is $\ZZ$; in fact our constructions are economical, in the sense of the following theorem, which is the subject of Section~\ref{section.construction}:

\begin{theorem}\label{theorem.construction}
For every integer $x$, there is a set $S\subseteq \{0,1,\dots,17|x|\}$ with $|S+S|-|S-S|=x$.
\end{theorem}

{\it Acknowledgements.} The first author was supported in part by grants from the Natural Sciences
and Engineering Research Council. The second author was supported in part by a grant from The
City University of New York PSC-CUNY Research Award Program. The second author also acknowledges
helpful discussions with Natella V. O'Bryant.

\section{Sums and differences in randomly chosen sets}

In this section, we establish several ancillary results on the probabilities that particular sums
and differences are present or absent in sets chosen randomly from certain classes of sets. We will
consider in particular the following classes: Let $n$, $\ell$, and $u$ be integers with
$n\ge\ell+u$. Fix $L\subseteq \{0,\dots,\ell-1\}$ and $U\subseteq \{n-u,\dots,n-1\}$. We will
consider the set of all subsets $A\subseteq \{0,\dots,n-1\}$ satisfying $A\cap \{0,\dots,\ell-1\} =
L$ and $A\cap \{n-u,\dots,n-1\} = U$ as a probability space endowed with the uniform probability,
where each such set $A$ occurs with the probability $2^{-(n-\ell-u)}$.

All of the calculations in this section are straightforward, but the details depend upon the size
and sometimes the parity of the particular sum or difference we are investigating, and so the
lemmas herein are rather ugly. The reader with limited tolerance could scan
Propositions~\ref{prop.sums} and~\ref{prop.differences} and move on to the next section without
significantly interrupting the flow of ideas.

We begin with three lemmas describing the probabilities of particular sums missing from $A+A$, where $A$ is chosen randomly from a class of the type indicated above.

\begin{lemma}
\label{lemma.one.small.sum}
Let $n$, $\ell$, and $u$ be integers with $n\ge\ell+u$.
Fix $L\subseteq \{0,\dots,\ell-1\}$ and $U\subseteq \{n-u,\dots,n-1\}$.
Suppose that $R$ is a uniformly randomly chosen subset of $\{\ell,\dots,n-u-1\}$, and set $A:=L\cup R\cup U$. Then for any integer $k$ satisfying $2\ell-1\le k\le n-u-1$, the probability
    \[
    \Prob{k\notin A+A} = \begin{cases}
    \big(\tfrac12\big)^{|L|} \big(\tfrac34\big)^{(k+1)/2-\ell}, & \text{if $k$ is odd,} \\
    \big(\tfrac12\big)^{|L|+1} \big(\tfrac34\big)^{k/2-\ell}, & \text{if $k$ is even.} \\
    \end{cases}
    \]
\end{lemma}

\begin{proof}
Define random variables $X_j$ by setting $X_j=1$ if $j\in A$ and $X_j=0$ otherwise. By the definition of $A$, the variables $X_j$ are independent random variables for $\ell\le j\le n-u-1$, each taking the values 0 and 1 with probability $1/2$ each, while the variables $X_j$ for $0\le j\le\ell-1$ and $n-u\le j\le n-1$ have values that are fixed by the choices of $L$ and $U$.

We have $k\notin A+A$ if and only if $X_j X_{k-j} = 0$ for all $0 \le j \le k/2$; the key point is that these variables $X_j X_{k-j}$ are independent of one another. Therefore
    \[
    \Prob{k\notin A+A} = \prod_{0\le j\le k/2} \Prob{X_j X_{k-j} = 0}.
    \]
If $k$ is odd, this becomes
    \begin{align*}
    \Prob{k\notin A+A} &= \prod_{j=0}^{\ell-1} \Prob{X_j X_{k-j} = 0} \prod_{j=\ell}^{(k-1)/2} \Prob{X_j X_{k-j} = 0} \\
    &= \prod_{j\in L} \Prob{X_{k-j} = 0} \prod_{j=\ell}^{(k-1)/2} \Prob{X_j = 0 \text{ or } X_{k-j} = 0} \\
    &= \big(\tfrac12\big)^{|L|} \big(\tfrac34\big)^{(k+1)/2-\ell}.
    \end{align*}
On the other hand, if $k$ is even then
    \begin{align*}
    \Prob{k\notin A+A} &= \prod_{j=0}^{\ell-1} \Prob{X_j X_{k-j} = 0} \bigg( \prod_{j=\ell}^{k/2-1} \Prob{X_j X_{k-j} = 0} \bigg) \Prob{X_{k/2}X_{k/2} = 0} \\
    &= \prod_{j\in L} \Prob{X_{k-j} = 0} \bigg( \prod_{j=\ell}^{k/2-1} \Prob{X_j = 0 \text{ or } X_{k-j} = 0} \bigg) \Prob{X_{k/2} = 0} \\
    &= \big(\tfrac12\big)^{|L|} \big(\tfrac34\big)^{k/2-\ell} \cdot \tfrac12.
    \end{align*}
\end{proof}

\begin{lemma}
\label{lemma.one.big.sum}
Let $n$, $\ell$, and $u$ be integers with $n\ge\ell+u$.
Fix $L\subseteq \{0,\dots,\ell-1\}$ and $U\subseteq \{n-u,\dots,n-1\}$.
Suppose that $R$ is a uniformly randomly chosen subset of $\{\ell,\dots,n-u-1\}$, and set $A:=L\cup R\cup U$. Then for any integer $k$ satisfying $n+\ell-1\le k\le 2n-2u-1$, the probability
    \[
    \Prob{k\notin A+A} = \begin{cases}
    \big(\tfrac12\big)^{|U|} \big(\tfrac34\big)^{n-(k+1)/2-u}, & \text{if $k$ is odd,} \\
    \big(\tfrac12\big)^{|U|+1} \big(\tfrac34\big)^{n-1-k/2-u}, & \text{if $k$ is even.} \\
    \end{cases}
    \]
\end{lemma}

\begin{proof}
This follows from Lemma \ref{lemma.one.small.sum} applied to the parameters $\ell'=u$ and $L'=n-1-U$, $u'=\ell$ and $U'=n-1-L$, and $A'=n-1-A$ and $k'=2n-2-k$.
\end{proof}

\begin{lemma}
\label{lemma.one.sum.totally.random}
Suppose that $A$ is a uniformly randomly chosen subset of $\{0,\dots,n-1\}$. Then for any integer $0\le k\le n-1$, the probability
    \[
    \Prob{k\notin A+A} = \begin{cases}
    \big(\tfrac34\big)^{(k+1)/2}, & \text{if $k$ is odd,} \\
    \tfrac12 \big(\tfrac34\big)^{k/2}, & \text{if $k$ is even;} \\
    \end{cases}
    \]
while for any integer $n-1\le k\le 2n-2$, the probability
    \[
    \Prob{k\notin A+A} = \begin{cases}
    \big(\tfrac34\big)^{n-(k+1)/2}, & \text{if $k$ is odd,} \\
    \tfrac12 \big(\tfrac34\big)^{n-1-k/2}, & \text{if $k$ is even.} \\
    \end{cases}
    \]
\end{lemma}

\begin{proof}
This follows immediately from Lemmas \ref{lemma.one.small.sum} and \ref{lemma.one.big.sum} upon setting $\ell=u=0$ and $L=U=\emptyset$.
\end{proof}

We now use these lemmas to establish the following proposition, in which we want a positive probability that many integers $k$ appear in the sumset $A+A$. While these events, varying over $k$, are not independent, we need only a lower bound on the probability; hence it suffices to combine crudely the exact probabilities given in Lemmas \ref{lemma.one.small.sum} and \ref{lemma.one.big.sum}. We emphasize that we have made no effort to optimize the lower bound given in the following proposition.

\begin{prop}
\label{prop.sums}
Let $n$, $\ell$, and $u$ be integers with $n\ge\ell+u$.
Fix $L\subseteq \{0,\dots,\ell-1\}$ and $U\subseteq \{n-u,\dots,n-1\}$.
Suppose that $R$ is a uniformly randomly chosen subset of $\{\ell,\dots,n-u-1\}$, and set $A:=L\cup R\cup U$. Then the probability that
    \[
    \{2\ell-1,\dots,n-u-1\} \cup \{n+\ell-1,\dots,2n-2u-1\} \subseteq A+A
    \]
is greater than $1-6(2^{-|L|}+2^{-|U|})$.
\end{prop}

\begin{proof}
We employ the crude inequality
    \begin{multline*}
    \Prob{\{2\ell-1,\dots,n-u-1\} \cup \{n+\ell-1,\dots,2n-2u-1\} \not\subseteq A+A} \\
        \le \sum_{k=2\ell-1}^{n-u-1} \Prob{k\notin A+A} + \sum_{k=n+\ell-1}^{2n-2u-1} \Prob{k\notin A+A}.
    \end{multline*}
The first sum can be bounded, using Lemma \ref{lemma.one.small.sum}, by
    \begin{align*}
    \sum_{k=2\ell-1}^{n-u-1} \Prob{k\notin A+A} &< \sum_{\substack{k\ge2\ell-1 \\ k\text{ odd}}} \big(\tfrac12\big)^{|L|} \big(\tfrac34\big)^{(k+1)/2-\ell} + \sum_{\substack{k\ge2\ell-1 \\ k\text{ even}}} \big(\tfrac12\big)^{|L|+1} \big(\tfrac34\big)^{k/2-\ell} \\
    &= \big(\tfrac12\big)^{|L|} \sum_{m=0}^\infty \big(\tfrac34\big)^m + \big(\tfrac12\big)^{|L|+1} \sum_{m=0}^\infty \big(\tfrac34\big)^m = 6\big(\tfrac12\big)^{|L|}.
    \end{align*}
The second sum can be bounded in a similar way using Lemma \ref{lemma.one.big.sum}, yielding
    \[
    \sum_{k=n+\ell-1}^{2n-2u-1} \Prob{k\notin A+A} < 6\big(\tfrac12\big)^{|U|}.
    \]
Therefore $\Prob{\{2\ell-1,\dots,n-u-1\} \cup \{n+\ell-1,\dots,2n-2u-1\} \not\subseteq A+A}$ is bounded above by $6(1/2)^{|L|} + 6(1/2)^{|U|}$, which is equivalent to the statement of the proposition.
\end{proof}

We turn now to three lemmas describing the probabilities that particular differences are missing
from $A-A$, where $A$ is chosen randomly from one of our classes. A new obstacle appears: while the
random variables $X_jX_{k-j}$ controlling the presence of the sum $k$ in $A+A$ are always mutually
independent, the same is not true of the random variables $X_jX_{k+j}$ controlling the presence of
the difference $k$ in $A-A$, at least when $k$ is small enough that $j$, $k+j$, and $2k+j$ can all
lie between 0 and $n-1$. Fortunately, when $k$ is this small the probabilities in question are
already minuscule, so a simple argument provides a serviceable bound
(Lemma~\ref{lemma.one.small.diff} below).

\begin{lemma}
\label{lemma.one.big.diff}
Let $n$, $\ell$, and $u$ be integers with $n\ge\ell+u$.
Fix $L\subseteq \{0,\dots,\ell-1\}$ and $U\subseteq \{n-u,\dots,n-1\}$.
Suppose that $R$ is a uniformly randomly chosen subset of $\{\ell,\dots,n-u-1\}$, and set $A:=L\cup R\cup U$. Then for any integer $k$ satisfying $n/2\le k\le n-u-\ell$, the probability
    \[
    \Prob{k\notin A-A} = \big(\tfrac12\big)^{|L|+|U|} \big(\tfrac34\big)^{n-\ell-u-k}.
    \]
\end{lemma}

\begin{proof}
Define random variables $X_j$ by setting $X_j=1$ if $j\in A$ and $X_j=0$ otherwise, as in the proof of Lemma \ref{lemma.one.small.sum}. We have $k\notin A-A$ if and only if $X_j X_{k+j} = 0$ for all $0 \le j \le n-1-k$, and again these variables $X_j X_{k+j}$ are independent of one another. Therefore
    \begin{align*}
    \Prob{k\notin A-A} &= \prod_{j=0}^{n-1-k} \Prob{X_j X_{k+j} = 0} \\
    &= \prod_{j=0}^{\ell-1} \Prob{X_j X_{k+j} = 0} \prod_{j=\ell}^{n-u-1-k} \Prob{X_j X_{k+j} = 0} \prod_{j=n-u-k}^{n-1-k} \Prob{X_j X_{k+j} = 0}  \\
    &= \prod_{j\in L} \Prob{X_{k+j} = 0} \prod_{j=\ell}^{n-u-1-k} \Prob{X_j = 0 \text{ or } X_{k+j} = 0} \prod_{j\in U-k} \Prob{X_j = 0} \\
    &= \big(\tfrac12\big)^{|L|} \big(\tfrac34\big)^{n-\ell-u-k} \big(\tfrac12\big)^{|U|}.
    \end{align*}
\end{proof}

\begin{lemma}
\label{lemma.one.small.diff}
Let $a$ and $b$ be integers with $a<b$.
Suppose that $R$ is a uniformly randomly chosen subset of $\{a,\dots,b-1\}$. Then for any integer $k$ satisfying $1\le k\le 2(b-a)/3$, the probability
    \[
    \Prob{k\notin R-R} \le \big(\tfrac34\big)^{(b-a)/3}.
    \]
\end{lemma}

\begin{remark}
In fact, the probability in question can be written exactly in terms of products of Fibonacci numbers: in the simplest case, $\Prob{1\notin R-R} = F_{b-a+2}/2^{b-a}$. However, the resulting expressions would become too tedious to handle in our applications below. When $b-a$ is large and $k$ is small, the actual value of the probability $\Prob{k\notin R-R}$ is proportional to $((1+\sqrt5)/4)^{b-a-k} \approx 0.809^{b-a}$, whereas the bound in Lemma~\ref{lemma.one.big.diff} gives $(3/4)^{(b-a)/3}\approx0.909^{b-a}$. However, in the particular case $k=(b-a)/2$, the probability in question is exactly $(3/4)^{(b-a)/2}\approx0.866^{b-a}$, so the bound in Lemma~\ref{lemma.one.big.diff} is not too unreasonable.
\end{remark}

\begin{proof}
Define the set
    \[
    J := \big\{ a\le j< b-k : \big\lfloor \tfrac{j-a}k \big\rfloor \text{ is even} \big\}.
    \]
In other words, $J$ contains the first $k$ integers starting at $a$, then omits the following $k$ integers, then contains the next $k$ integers, and so on until the upper bound $a+2(b-a)/3$ is reached. The following properties of $J$ can be easily verified:
\begin{enumerate}
\renewcommand{\labelenumi}{(\roman{enumi})}
\item if $j\in J$, then $j+k\notin J$;
\item $|J| \ge (b-a)/3$.
\end{enumerate}
Now define random variables $X_j$ by setting $X_j=1$ if $j\in R$ and $X_j=0$ otherwise, as in the
proof of Lemma~\ref{lemma.one.small.diff}. We have $k\notin R-R$ if and only if $X_j X_{k+j} = 0$
for all $a \le j < b-k$.
    \begin{align*}
    \Prob{k\notin R-R} &= \Prob{X_j X_{k+j} = 0 \text{ for all }a \le j < b-k} \\
    &\le \Prob{X_j X_{k+j} = 0 \text{ for all }j\in J}.
    \end{align*}
However, property~(i) above ensures that the random variables $X_j X_{k+j}$ are independent of one
another as $j$ ranges over $J$. Therefore
    \[
    \Prob{k\notin R-R}      \le \prod_{j\in J} \Prob{X_j X_{k+j} = 0}
                            = \big( \tfrac34 \big) ^{|J|} \le \big( \tfrac34 \big) ^{(b-a)/3}
    \]
by property (ii) above.
\end{proof}

\begin{lemma}
\label{lemma.one.diff.totally.random}
Suppose that $A$ is a uniformly randomly chosen subset of $\{0,\dots,n-1\}$. Then for any integer $1\le k\le n/2$, the probability $\Prob{k\notin A-A} \le (3/4)^{n/3}$, while for any integer $n/2\le k\le n-1$, the probability $\Prob{k\notin A-A} = (3/4)^{n-k}$.
\end{lemma}

\begin{proof}
The first assertion follows immediately from Lemma~\ref{lemma.one.small.diff} upon setting $a=0$ and $b=n$, while the second assertion follows immediately from Lemma~\ref{lemma.one.big.diff} upon setting $\ell=u=0$ and $L=U=\emptyset$.
\end{proof}

We now use these lemmas to establish the following proposition, in which we want a positive
probability that many integers $k$ appear in the difference set $A-A$. Again it suffices to combine
crudely the results of Lemmas \ref{lemma.one.big.diff} and \ref{lemma.one.small.diff}, since we
need only a lower bound on the probability. Once again we have emphasized ease of exposition over
optimization of the lower bound itself; in particular, we could have achieved better constants at
the expense of uglier technicalities.

\begin{prop}
\label{prop.differences}
Let $n$, $\ell$, and $u$ be integers with $n\ge4(\ell+u)$.
Fix $L\subseteq \{0,\dots,\ell-1\}$ and $U\subseteq \{n-u,\dots,n-1\}$.
Suppose that $R$ is a uniformly randomly chosen subset of $\{\ell,\dots,n-u-1\}$, and set $A:=L\cup R\cup U$. Then the probability that
   \[
    \{-(n-\ell-u),\dots,n-\ell-u\} \subseteq A-A
    \]
is greater than $1- 4(1/2)^{|L|+|U|}-(n/2)(3/4)^{(n-\ell-u)/3}$.
\end{prop}

\begin{proof}
By the symmetry of $A-A$ about 0 and the fact that $0\in A-A$ for any nonempty set $A$, it suffices to show that $A-A$ contains $\{1,\dots,n-\ell-u\}$. We employ the crude inequality
    \begin{align*}
    \Prob{\{1,\dots,n-\ell-u\} \not\subseteq A-A} &\le \sum_{k=1}^{n-\ell-u} \Prob{k\notin A-A} \\
    &\le \sum_{1\le k\le n/2} \Prob{k\notin R-R} + \sum_{n/2<k\le n-\ell-u} \Prob{k\notin A-A}.
    \end{align*}
The first sum can be bounded using Lemma \ref{lemma.one.small.diff} with $a=\ell$ and $b=n-u$; it is here that we use the hypothesis $n\ge4(\ell+u)$, to guarantee that every $k$ in the range $1\le k\le n/2$ satisfies $k\le2(n-\ell-u)/3$. We obtain
    \[
    \sum_{1\le k\le n/2} \Prob{k\notin R-R} \le \tfrac n2 \big( \tfrac34 \big)^{(n-\ell-u)/3}.
    \]
The second sum can be bounded using Lemma \ref{lemma.one.big.sum}, yielding
    \[
    \sum_{n/2<k\le n-\ell-u} \Prob{k\notin A-A} < \sum_{k=-\infty}^{n-\ell-u} \big(\tfrac12\big)^{|L|+|U|} \big(\tfrac34\big)^{n-\ell-u-k} = 4\big(\tfrac12\big)^{|L|+|U|}.
    \]
Therefore $\Prob{\{-(n-\ell-u),\dots,n-\ell-u\} \not\subseteq A-A}$ is bounded above by $4(1/2)^{|L|+|U|} + (n/2)(3/4)^{(n-\ell-u)/3}$, which is equivalent to the statement of the proposition.
\end{proof}


\section{Proof of Theorem \ref{positive.proportion.thm}}

In this section we show that the collections of \sd\ sets, \dd\ sets, and \sdb\ sets all have positive lower density. Our strategy is to fix the ``fringes'' of a subset of $\nset$ (that is, stipulate which integers close to 0 and $n-1$ are and are not in the set) in a way that forces the set to have missing differences (or sums). We then use the probabilistic lemmas of the previous section to show that for many sets with the prescribed fringes, all other sums (or differences) will be present. We have not attempted to optimize the constants appearing in the following three theorems, in part because the previous section would have become even more technical and ugly, and in part because we were unlikely to have come close to the true constants (see Conjecture~\ref{three.rho.values.conjecture} below) in any event.

We begin by showing that a positive proportion of sets are \sd. Here, choosing appropriate fringes is most non-trivial, compared to the two theorems that follow.

\newcommand{\sdconst}{(2\times 10^{-7})}
\begin{theorem}
\label{many.sd.sets.theorem}
For $n\ge15$, the number of \sd\ subsets of $\nset$ is at least $\sdconst 2^n$.
\end{theorem}

\begin{proof}
First, note that the bound $\sdconst 2^n$ is less than 1 for $15\le n\le22$; the existence of the single \sd\ set $\{0,2,3,4,7,11,12,14\}$ is enough to verify the theorem in that range. Henceforth we can assume that $n\ge23$.

Define $L:=\{0,2,3,7,8,9,10\}$ and $U:=\{n-11,n-10,n-9,n-8,n-6,n-3,n-2,n-1\}$. We show that the
number of \sd\ subsets $A\subseteq\nset$ satisfying $A\cap\{0,\dots,10\} = L$ and
$A\cap\{n-11,\dots,n-1\}=U$ is at least $\sdconst 2^n$. For any such $A$, the fact that $U-L$ does
not contain $n-7$ implies that $A-A$ contains neither $n-7$ nor $-(n-7)$; since
$A-A\subseteq\{-(n-1),\dots,n-1\}$, we see that
    \[|A-A|\le 2n-3.\]
Therefore it suffices to show that there are at least $\sdconst 2^n$ sets $A$, satisfying $A\cap\{0,\dots,10\} = L$ and $A\cap\{n-11,\dots,n-1\}=U$, for which $|A+A|\ge2n-2$.

For any such $A$, we see by direct calculation that $A+A$ contains
    \begin{align*}
    L+L&= \{0,\dots,20\} \setminus \{1\}, \\
    L+U&=\{n-11,\dots,n+9\}, \\
    U+U&=\{2n-22,\dots,2n-2\}.
    \end{align*}
In particular, if $23\le n\le32$ then $A+A$ automatically equals $\{0,\dots,2n-2\} \setminus \{1\}$, giving $|A+A|=2n-2$; the number of such $A$ is exactly $2^{n-22} > \sdconst 2^n$, since there are $n-22$ numbers between 11 and $n-12$ inclusive.

For $n\ge33$, Proposition~\ref{prop.sums} (applied with $\ell=r=11$) tells us that when $A$ is chosen uniformly randomly from all such sets, the probability that $A+A$ contains $\{21,\dots,n-12\} \cup \{n+10,\dots,2n-23\}$ is at least
\[
1-6(2^{-|L|}+2^{-|U|}) = 1-6(2^{-7}+2^{-8}) = \frac{119}{128}.
\]
In other words, there are at least $2^{n-22}\cdot119/128 > \sdconst 2^n$ such sets $A$. For all these sets, we see that $A+A$ again equals $\{0,\dots,2n-2\} \setminus \{1\}$, and hence all such sets are \sd.
\end{proof}

The next two theorems carry out a similar approach to showing that a positive proportion of sets are \dd\ or \sdb. These two results appeal to the serviceable but crude Lemma \ref{lemma.one.small.diff}, and consequently the constants that appear, as well as the computation needed to take care of smaller values of $n$, are likewise far from optimal.

\newcommand{\ddconst}{0.0015\cdot}
\begin{theorem}
\label{many.dd.sets.theorem}
For $n\ge4$, the number of \dd\ subsets of $\nset$ is at least $\ddconst 2^n$.
\end{theorem}

\begin{proof}
The bound can be verified computationally for small $n$: we have computed by exhaustive search for $n\le 27$ the
number of \dd\ subsets $\nset$ that contain both 0 and $n-1$. Counting just these sets and their
translates is enough to prove this theorem for $n\le 39$. Henceforth, we assume that $n\ge 40$.

Define $L:=\{0,2,3\}$ and $U:=\{n-2,n-1\}$. We show that the number of \dd\ subsets
$A\subseteq\nset$ satisfying $A\cap\{0,1,2,3\} = L$ and $A\cap\{n-2,n-1\}=U$ is at least $\ddconst
2^n$. For any such $A$, the fact that $L+L$ does not contain 1 implies that $A+A$ does not contain
1, and so $|A+A|\le 2n-2$. Therefore it suffices to show that there are at least $\ddconst 2^n$
sets $A$, satisfying $A\cap\{0,1,2,3\} = L$ and $A\cap\{n-2,n-1\}=U$, for which $|A-A|=2n-1$.

For any such $A$, we see by direct calculation that $A-A$ contains
    \[
    (L-U) \cup (U-L) = \{-(n-5),\dots,-(n-1)\}\cup\{n-5,\dots,n-1\}.
    \]
Furthermore, Proposition~\ref{prop.differences} (applied with $\ell=4$, $u=2$, and $n\ge 24$) tells
us that when $A$ is chosen uniformly randomly from all such sets, the probability that $A-A$
contains $\{-(n-6),\dots,n-6\}$ is at least
\begin{align*}
1- 4\bigg( \frac12 \bigg) ^{|L|+|U|}-\bigg( \frac n2 \bigg) \bigg( \frac34 \bigg)^{(n-\ell-u)/3} &= 1- 4\bigg( \frac12 \bigg)^5-\bigg( \frac n2 \bigg) \bigg( \frac34 \bigg)^{(n-6)/3} \\
&= \frac78 - \frac{8n}9\bigg( \frac34 \bigg)^{n/3}.
\end{align*}
As a function of $n$, this expression is increasing for $n\ge11$, and at $n=40$ its value is larger
than $0.107536$. In other words, there are at least $2^{n-6}\cdot 0.107536 > \ddconst 2^n$ such
sets $A$. For all these sets, we see that $A-A$ equals $\{-(n-1),\dots,n-1\}$, and hence all such
sets are \dd.
\end{proof}

\newcommand{\sdbconst}{(2\times 10^{-5})}
\begin{theorem}
\label{many.sdb.sets.theorem}
For $n\ge1$, the number of \sdb\ subsets of $\nset$ is at least $\sdbconst 2^n$.
\end{theorem}

\begin{proof}
The bound can be verified computationally for small $n$: for $n\le 27$ we have computed the exact number of \sdb\
subsets of $\nset$ that contain both $0$ and $n-1$. Counting only these sets and their translates
proves the theorem for $n\le 42$. Henceforth, we assume that $n\ge43$.

Define $L:=\{0,\dots,5\}$ and $U:=\{n-6,\dots,n-1\}$. We give a lower bound for the number of \sdb\
subsets $A\subseteq\nset$ satisfying $L\cup U\subseteq A$; in fact we show that the number of such
subsets with $|A+A|=|A-A|=2n-1$, the maximum possible size, is at least $\sdbconst 2^n$. Combining
Propositions~\ref{prop.sums} and~\ref{prop.differences} (applied with $\ell=u=6$), we find that
when $A$ is chosen uniformly randomly from all such sets, the probability that both $A+A$ and $A-A$
are as large as possible is at least
    \[
    1- 6(2^{-|L|} + 2^{-|R|}) - 4\bigg( \frac12 \bigg) ^{|L|+|U|}
        -\bigg( \frac n2 \bigg) \bigg( \frac34 \bigg)^{(n-\ell-u)/3}
        = \frac34 - \frac{8n}9\bigg( \frac34 \bigg)^{n/3}.
    \]
This function is increasing for $n\ge1$ and takes a value larger than 0.131232 when $n=43$. In
other words, there are at least $2^{n-12}\cdot 0.131232 > \sdbconst 2^n$ such sets $A$. For all
these sets, we see that $A+A$ equals $\{0,\dots,2n-2\}$ and $A-A$ equals $\{-(n-1),\dots,n-1\}$,
and hence all such sets are \sdb.
\end{proof}

\section{Average values}\label{section.average.values}

In this section, we prove Theorem~\ref{theorem.average.value} by calculating the average values of
$|S-S|$ and $|S+S|$ as $S$ ranges over an arithmetic progression $P$ of length $n$. Since the
problem is invariant under dilations and translations, it suffices to prove the theorem in the case
$P=\nset$.

We begin by addressing the average cardinality of the sumset $S+S$. In fact, we can give an exact
formula for the average size of the sumset, or equivalently for the sum of the sizes of all sumsets
as $S$ ranges over subsets of $\nset$. The reason we can do so is essentially because of the
linearity of expectations of random variables.

\begin{theorem}
For any positive integer $n$, we have
\begin{equation}
\label{sum.average.equation}
\sum_{S\subseteq\smallnset} |S+S| = 2^n(2n-11) + \begin{cases}
19\cdot3^{(n-1)/2}, &\text{if $n$ is odd,} \\
11\cdot3^{n/2}, &\text{if $n$ is even.}
\end{cases}
\end{equation}
\label{sum.average.theorem}
\end{theorem}

\begin{proof}
We begin with the manipulation
\begin{multline}
\sum_ {S\subseteq \smallnset} |S+S|
    = \sum_ {S\subseteq \smallnset} \sum_{\substack{0\le k\le2n-2 \\ k\in S+S}} 1
    = \sum_{k=0}^{2n-2} \sum_{\substack{S\subseteq \smallnset \\ k\in S+S}} 1 \\
    = \sum_{k=0}^{2n-2} 2^n \Prob{k\in S+S}
    = 2^n(2n-1) - 2^n \sum_{k=0}^{2n-2}\Prob{k\notin S+S}.
\label{before.even.odd.equation}
\end{multline}
We suppose that $n=2m+1$ is odd, the case where $n$ is even being similar. We begin by considering only the lower half of possible values for $k$. By Lemma~\ref{lemma.one.sum.totally.random}, we have
    \begin{align*}
\sum_{k=0}^{n-2} \Prob{k\notin S+S}
    &= \sum_{j=0}^{m-1} \Prob{2j\notin S+S} + \sum_{j=0}^{m-1} \Prob{2j+1\notin S+S} \\
    &= \sum_{j=0}^{m-1} \tfrac12 \big( \tfrac34 \big)^j + \sum_{j=0}^{m-1} \big( \tfrac34 \big)^{j+1} = 5\big( 1- \big( \tfrac34 \big)^m \big).
    \end{align*}
By the symmetry of $S+S$ about $n-1$, the same calculation holds for $\sum_{k=n}^{2n-2} \Prob{k\notin S+S}$. Therefore, appealing to Lemma~\ref{lemma.one.sum.totally.random} again for $k=n-1=2m$,
\[
\sum_{k=0}^{2n-2}\Prob{k\notin S+S} = 5\big( 1- \big( \tfrac34 \big)^m \big) + \tfrac12 \big( \tfrac34 \big)^m + 5\big( 1- \big( \tfrac34 \big)^m \big) = 10 - \tfrac{19}2 \big( \tfrac34 \big)^{(n-1)/2}.
\]
Inserting this value into the right-hand side of equation~\eqref{before.even.odd.equation} establishes the lemma for odd~$n$. A similar calculation gives the result for even~$n$.
\end{proof}

While it is possible to write down an exact formula for the average size of the difference set $S-S$ as $S$ ranges over all subsets of $\nset$, the formula would be far too ugly to be of use. We prefer in this case to present a simple asymptotic formula with a reasonable error term.

\begin{theorem}
For any positive integer $n$, we have
\begin{equation}
\label{difference.average.equation}
\sum_{S\subseteq\smallnset} |S-S| = 2^n(2n-7) + O\big( n6^{n/3} \big).
\end{equation}
\label{difference.average.theorem}
\end{theorem}

\begin{proof}
As in the proof of the previous theorem, we have
\begin{align}
\sum_ {S\subseteq \smallnset} |S-S|
   & = \sum_ {S\subseteq \smallnset} \sum_{\substack{-(n-1)\le k\le n-1 \\ k\in S-S}} 1
    = \sum_{k=-(n-1)}^{n-1} \sum_{\substack{S\subseteq \smallnset \\ k\in S-S}} 1 \notag \\
    &= \sum_{k=-(n-1)}^{n-1} 2^n \Prob{k\in S-S} \notag \\
    &= 2^n(2n-1) - 2^n \sum_{k=-(n-1)}^{n-1} \Prob{k\notin S-S} \notag \\
    &= 2^n(2n-1) -1- 2^{n+1} \sum_{k=1}^{n-1} \Prob{k\notin S-S},
\label{positive.k.equation}
\end{align}
the last equality following from the symmetry of $S-S$ around 0 and the fact that 0 is in $S-S$ for
nonempty $S$. From Lemma~\ref{lemma.one.diff.totally.random} we have
    \[
    \sum_{k=1}^{\lceil n/2\rceil-1} \Prob{k\notin S-S} \le \sum_{k=1}^ {\lceil n/2\rceil-1} \big( \tfrac34 \big)^{n/3} \le n \big( \tfrac34 \big)^{n/3}
    \]
and
    \[
    \sum_{k=\lceil n/2\rceil}^{n-1} \Prob{k\notin S-S} = \sum_{k=\lceil n/2\rceil}^{n-1} \big( \tfrac34 \big)^{n-k} = 3\big( 1-\big( \tfrac34 \big)^{n-\lceil n/2\rceil} \big),
    \]
which combine to give
    \[
    \sum_{k=1}^{n-1} \Prob{k\notin S-S} = 3 + O\big( n\big( \tfrac34 \big)^{n/3} \big).
    \]
Inserting this expression into the right-hand side of equation~\eqref{positive.k.equation} establishes the theorem.
\end{proof}

Examining the derivations of these two theorems reveals that it really is the commutativity $s_1+s_2=s_2+s_1$ that causes the difference in the average sizes of $S+S$ and $S-S$: a typical potential element of $S+S$ has only about half as many chances to be realized as a sum as the corresponding potential element of $S-S$ has at being realized as a difference. To further emphasize this observation, we note that if the single set $S$ is replaced by two sets $S$ and $T$, the disparity disappears: for an arithmetic progression $P$ of length $n$, we have
    \[
    \frac1{2^{2n}} \sum_{S\subseteq P} \sum_{T\subseteq P} |S-T| \sim \frac1{2^{2n}} \sum_{S\subseteq P} \sum_{T\subseteq P} |S+T| \sim 2n-7.
    \]

\section{Sets with prescribed imbalance between sums and differences}
\label{section.construction}

In this section we prove that the range of possible values for $|S+S|-|S-S|$ is all of $\Z$. Furthermore, as asserted in Theorem~\ref{theorem.construction}, our constructions show that for every integer $x$, we can find a subset $S$ of $\{0,\dots,17|x|\}$ such that $|S+S|-|S-S|=x$. As one might expect from the foregoing discussion, the case where $x$ is negative is easiest.

\medskip\noindent{\it Negative values of $x$.}
For any integer $x<0$, set $S_x=\{0,\dots,|x|+1\}\cup\{2|x|+2\}$. Then
$S_x+S_x=\{0,1,\dots,3|x|+3\}\cup\{4|x|+4\}$ while $S_x-S_x=\{-(2|x|+2),\dots,2|x|+2\}$, whereupon
\[
|S_x+S_x|-|S_x-S_x| = (3|x|+5) - (4|x|+5) = -|x| = x.
\]
Even more generally, take any integer $n\ge|x|+2$ and set $S=\{0,\dots,n-1\} \cup \{n+|x|\}$. Then $S+S=\{0,\dots,2n+|x|-1\} \cup \{2n+2|x|\}$ and $S-S=\{-(n+|x|), \dots, n+|x|\}$, which again yields $|S+S|-|S-S|=x$.
\medskip

We turn now to nonpositive values of $x$. Our general construction works for larger values of $x$, but we need to handle a few small values of $x$ individually.

\medskip\noindent{\it A few special cases.}
For a few small values of $x$, we find suitable sets $S_x$ simply by computation: if we set
\begin{align}
S_0 &:= \emptyset \notag \\
S_1 &:= \{0,2,3,4,7,11,12,14\} \label{S1.definition} \\
S_2 &:= \{0, 1, 2, 4, 5, 9, 12, 13, 14, 16, 17\} \notag \\
S_4 &:= \{0, 1, 2, 4, 5, 9, 12, 13, 17, 20, 21, 22, 24, 25\}, \notag
\end{align}
then in each case it can be checked that $|S_x+S_x|-|S_x-S_x| = x$. In fact, these examples are all
minimal in the sense that the diameter $\max S-\min S$ is as small as possible. (Vishaal Kapoor and
Erick Wong confirmed computationally the fact that $S_4$ is the unique, up to reflection, set of
integers of diameter at most 25 for which the sumset has four more elements than the difference
set. We note that Pigarev and Fre{\u\i}man~\cite{MR0434995} gave the slightly larger example
$S_4'=\{0,1,2,4,5,9,12,13,14,16,17,21,24,25,26,28,29\}$, which also satisfies $|S_4'+S_4'| -
|S_4'-S_4'| = 4$.)

In fact, these diameter-minimal examples are unique, up to reflection, except for $S_1$: there are
two other subsets of $\{0,\dots,14\}$, namely
\[
S_1' = \{0,1,2,4,5,9,12,13,14\}
\]
and its reflection, for which the sumset has one element more than the difference set. The first
set $S_1$ has only eight elements, as compared with the nine elements of $S_1'$. In fact,
Hegarty~\cite{hegarty} has shown that $S_1$ is also {\em the} \sd\ set with the smallest
cardinality, unique up to dilation, translation, and reflection. On the other hand, there are
tantalizing similarities among the sets $S_1'$, $S_2$, $S_4$, and $S_4'$ that might admit a clever
generalization.

We note that Ruzsa \cite{MR1171750} claimed that $U=\{0,1,3,4,5,6,7,10\}$ is \sd,
but this is incorrect: both $U+U$ and $U-U$ have 19 elements. We also mention the following observation of Hegarty: if one sets
\begin{align*}
A &= S_4 \cup (S_4 + 20) \\
&= \{0,1,2,4,5,9,12,13,17,20,21,22,24,25,29,32,33,37,40,41,42,44,45\},
\end{align*}
then one has $|A+A| = 91$ and $|A-A| = 83$, providing the statistic $\log 91 / \log 83 = 1.0208\dots$ which is important when using the elements of $A$ as ``digits''. More precisely, considering sets of the form $A_n = A+bA+b^2A+\dots+b^{n-1}A$ for suitably large fixed $b$, we have $|A_n+A_n| = |A_n-A_n|^{1.0208\dots}$, which is currently the best exponent known.

\medskip

For other positive values of $x$, the basic general construction is an adaptation of the subset
$S_1 \times \{0, \dots, k\}$ of $\Z\times\Z$, embedded in $\Z$ itself by a common technique of
regarding the coordinates as digits in a base-$b$ representation for suitably large $b$. In our
simple case, we can be completely explicit from the start.

\medskip\noindent{\it Odd values of $x$ exceeding 1.}
Let $x=2k+1$ with $k\ge1$. With $S_1$ defined as in equation~\eqref{S1.definition}, set
\begin{align}
S_{2k+1} &= S_1 + \{0, 29, 58, \dots, 29k \} \label{Sodd.definition} \\
&= \{ 0\le s\le 29k+14\colon s \equiv 0,2,3,4,5,11,12,\text{or }14 \mod{29} \}. \notag
\end{align}
Then we find that
\begin{align*}
S_{2k+1} + S_{2k+1} &= (S_1+S_1) + \{0, 29, 58, \dots, 58k\} \\
&= \{ 0\le s<29(2k+1) \colon s \not\equiv 1,20,\text{or }27 \mod{29} \},
\end{align*}
which reveals that $|S_{2k+1}+S_{2k+1}| = 26(2k+1)$. On the other hand,
\[
S_{2k+1} - S_{2k+1} = \big\{ {-29}\big(k+\tfrac12\big)<s<29\big(k+\tfrac12\big) \colon s \not\equiv -13,-6,6,\text{or }13 \mod{29} \big\},
\]
showing that $|S_{2k+1}-S_{2k+1}| = 25(2k+1)$, and so $|S_{2k+1}+S_{2k+1}|-|S_{2k+1}-S_{2k+1}| = 2k+1$ as desired.

\medskip\noindent{\it Even values of $x$ exceeding 4.}
Let $x=2k$ with $k\ge3$. With $S_{2k+1}$ defined as in equation~\eqref{Sodd.definition}, set $S_{2k} = S_{2k+1} \setminus \{29\}$. One can check that $S_{2k}-S_{2k}$ still equals all of $S_{2k+1}-S_{2k+1}$ but that $S_{2k}+S_{2k} = (S_{2k+1}+S_{2k+1}) \setminus \{29\}$. Therefore
\[
|S_{2k}+S_{2k}|-|S_{2k}-S_{2k}| = |S_{2k+1}+S_{2k+1}|-|S_{2k+1}-S_{2k+1}| -1 = 2k
\]
as desired. Notice that $S_{2k}$ is indeed contained in $\{0,\dots,17(2k)\}$ as asserted by Theorem~\ref{theorem.construction}, the closest call being the comparison between $\max S_6=101$ and $17\cdot 6=102$.


We note that as this manuscript was in preparation, Hegarty~\cite{hegarty}*{Theorem 9}
independently proved that $|S+S|-|S-S|$ can take all integer values $t$. In fact he proved,
extending ideas originating in our proof of Theorem~\ref{positive.proportion.thm}, somewhat more:
for each fixed integer $t$, if $n$ is sufficiently large then a positive proportion of subsets $S$
of $\nset$ satisfy $|S+S|-|S-S|=t$.

\section{Analysis of data}

\begin{figure}[hbtf]
\includegraphics[width=6in]{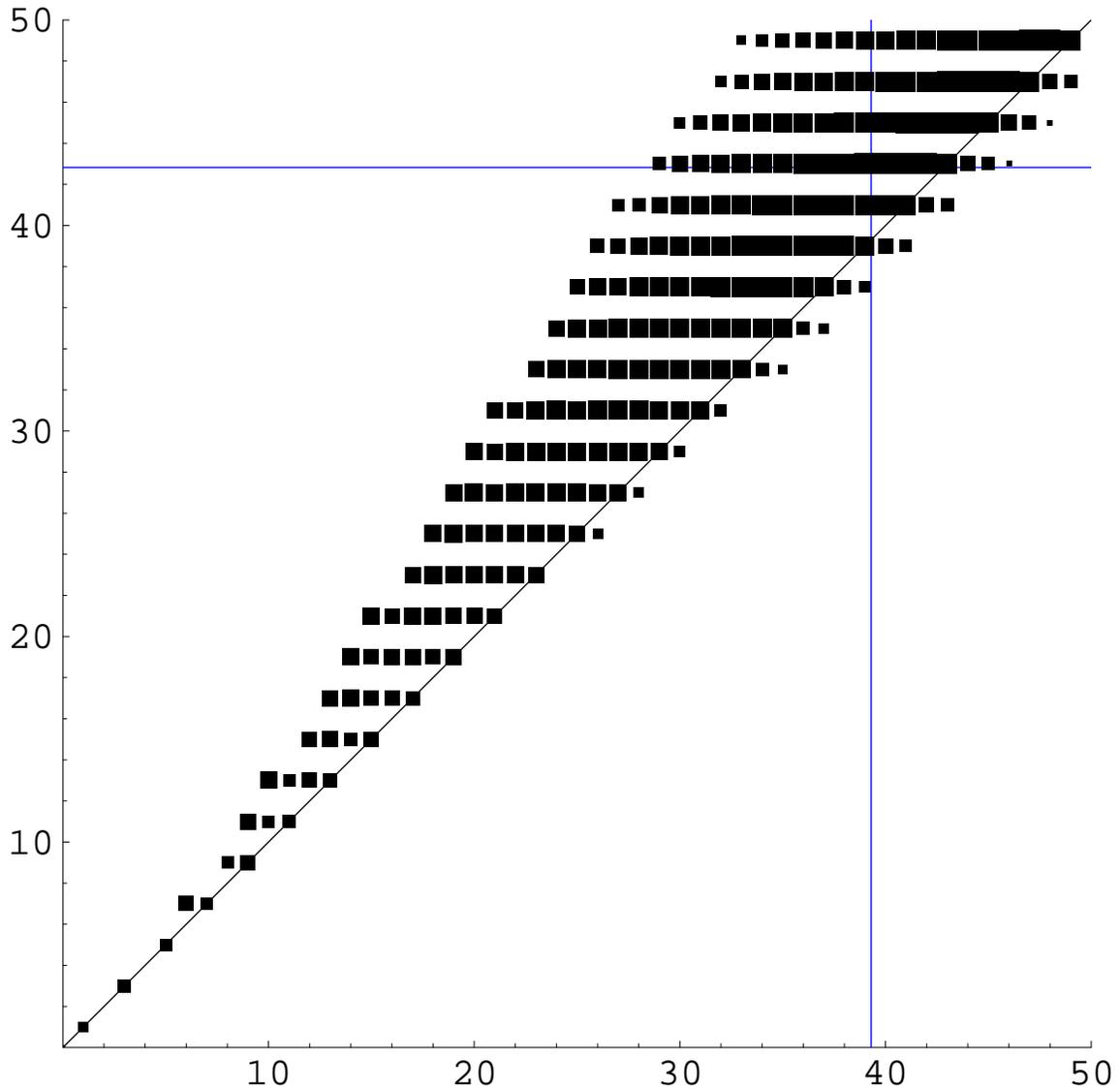}
\caption{The size of the square centered at $(x,y)$ indicates the number of subsets of
$\{0,\dots,24\}$ with $(|S+S|,|S-S|)=(x,y)$.}
\label{figure.n25}
\end{figure}

Theorem~\ref{theorem.average.value} gave the expected values of $|S+S|$ and $|S-S|$, which seems
most naturally phrased as saying that the expected number of missing sums is asymptotically 10,
while the expected number of missing differences is asymptotically~6. One is naturally led to
enquire as to the details of the joint distribution of these two quantities. Let $c_n(x,y)$ be the
number of subsets of $\nset$ with $|S+S|=x$ and $|S-S|=y$. Figure~\ref{figure.n25} shows a square
centered at $(x,y)\in\ZZ^2$ whose area is proportional to $\log(1+c_{25}(x,y))$. Also shown are the
lines $x=2^{-25}\sum_S |S+S|$ (the average size of a sumset), $y=2^{-25}\sum_S |S-S|$ (the average
size of a difference set), and $y=x$.

\begin{figure}[hbtf]
\includegraphics[width=6in]{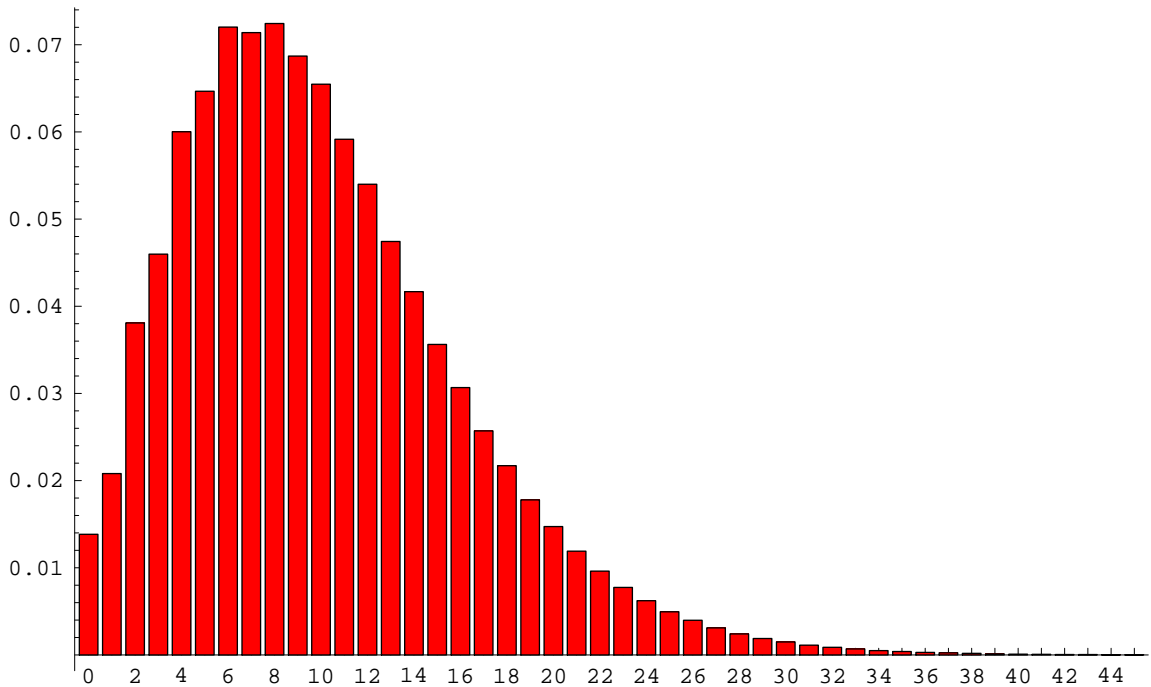}
\caption {The observed frequencies of the number of missing sums}
\label{missing.subsums}
\end{figure}

Figure~\ref {missing.subsums} shows the observed distribution of $X:=2n-1-|S+S|$ (that is, the
number of missing sums) for three million randomly generated subsets of $\{0,1,2,\dots,999\}$. For
example, the histogram shows that approximately 1.4\% of these subsets $S$ have the largest
possible sumset $S+S=\{0,\dots,1998\}$, approximately 2.1\% of them have exactly one element of
$\{0,\dots,1998\}$ missing from their sumsets, and so on. The histogram is essentially identical to
one generated from the complete data set for subsets of $\{0,\dots,26\}$.

Notice that there is a ``divot'' at the top of the histogram: the observed frequencies of sets
missing exactly 6 or 8 sums are both larger than the observed frequency of sets missing exactly 7
sums. In fact, the frequency for every even value seems to be larger than the average of its two
neighbors, while the opposite is true for the frequencies of the odd values; in other words, the
piecewise linear graph that connected the points at the tops of the histogram's bars would
alternate between being convex and concave.

Recall that the missing sums are typically very near the edges of the interval of possible sums. In
particular, the missing sums for a subset $S$ of $\{0,\dots,999\}$ tend to be near either 0 or
1998, and are therefore so far apart that their numbers are independent. Therefore the distribution
shown in Figure~\ref {missing.subsums} is the sum of two independent, identically distributed (by
symmetry) random variables that count the number of missing sums near one end. This is also
essentially the same distribution as the number $Y$ of missing sums in randomly chosen (infinite)
subsets $A$ of the nonnegative integers $\{0,1,\dots\}$. That is, if $Y_1,Y_2$ are independent with
the same distribution as $Y$, then $X$ and $Y_1+Y_2$ have approximately the same distribution (for
large $n$).

\begin{figure}[hbtf]
\includegraphics[width=3in]{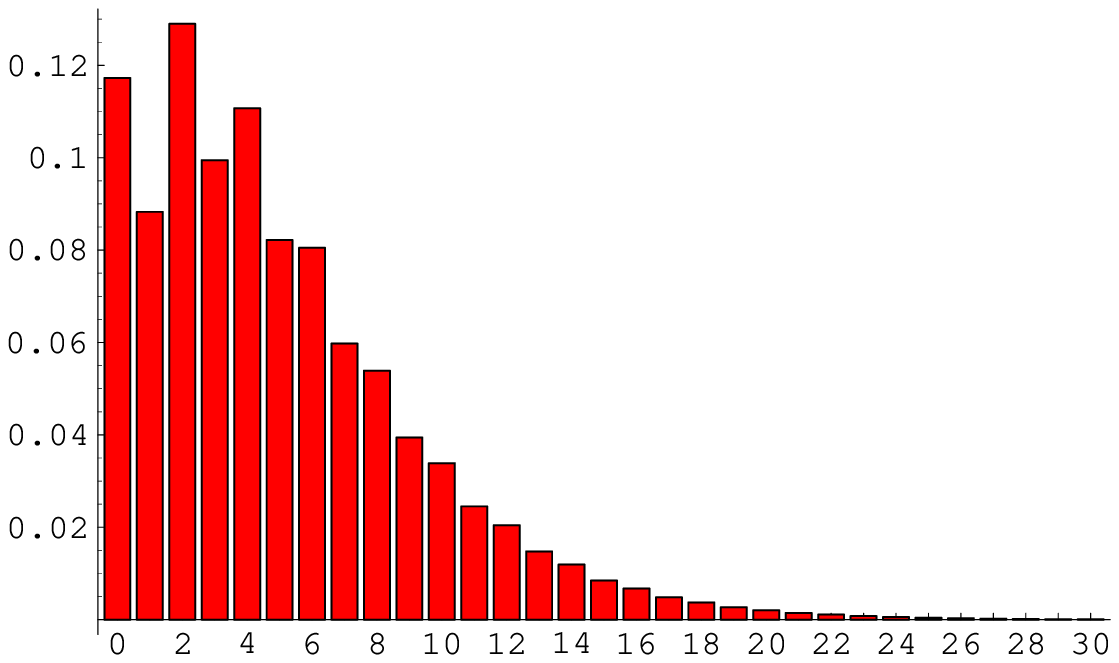}
\hfill
\includegraphics[width=3in]{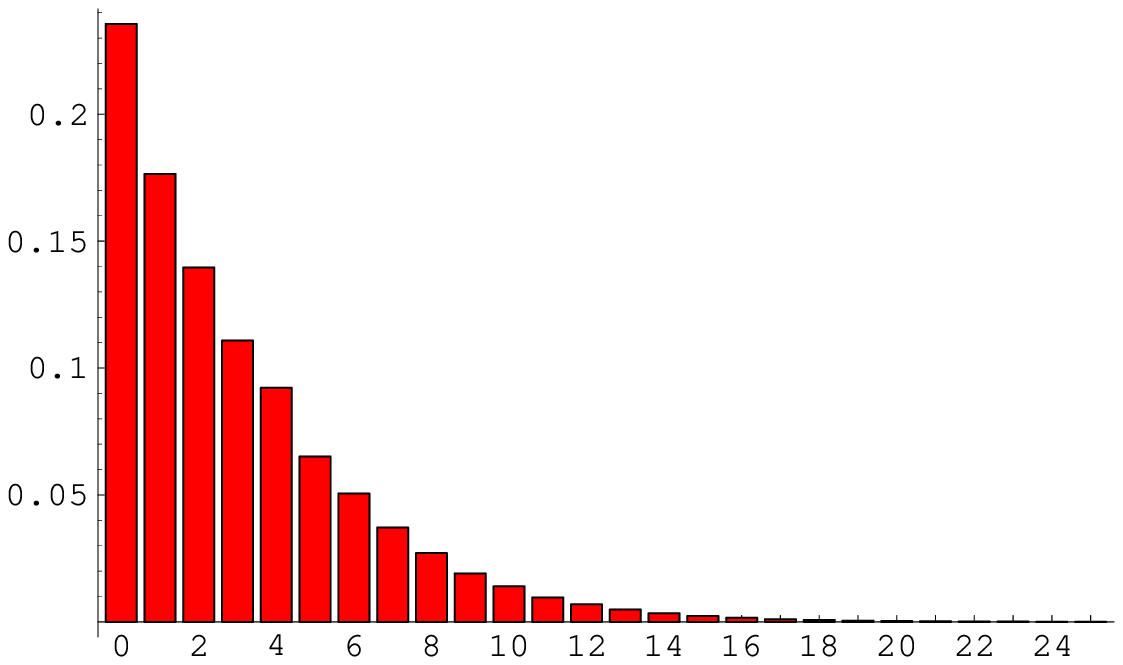}
\caption{The observed frequencies of the number of missing sums for randomly chosen subsets of $\{0,1,\dots\}$, with no restriction (left) and with the restriction that 0 belongs to the set (right)}
\label{one.sided.missing.subsums}
\end{figure}

At first one might think, then, that the parity phenomenon in Figure~\ref {missing.subsums} is
caused by that distribution being the sum of two independent copies of a simpler distribution.
However, in this latter distribution (the first histogram in Figure~\ref
{one.sided.missing.subsums}), the disparity between odd and even values is even more apparent.

Fortunately, the phenomenon here is easy to analyze: if 0 is not in our randomly chosen subset of
$\{0,1,\dots\}$, then there are automatically 2 missing sums, namely 0 and 1, and the rest of the
random subset can be shifted downwards by 1 to find the distribution of other missing sums:
    \[
    \Prob{Y=k}=\sum_{i=0}^{\lfloor k/2 \rfloor} \Prob{Y=k-2i \mid \min A=0} 2^{-i}.
    \]
In other words, there is a yet more fundamental distribution (the second histogram in Figure~\ref
{one.sided.missing.subsums}), given by the number of missing subsums in a randomly chosen subset of
$\{0,1,\dots\}$ containing~0. For example, that histogram shows that if a subset $S$ of
$\{0,1,\dots\}$ containing 0 is chosen at random, there is about a 23.6\% chance that
$S+S=\{0,1,\dots\}$.

The parity discrepancy seems to be absent in this last distribution, suggesting that it should be
the focus of further analysis; the two more complicated preceding distributions can be
reconstructed from suitable manipulations of this most fundamental one. The histogram suggests the
existence of a function $f(x)$, smooth and decaying faster than exponentially, such that the
probability of a randomly chosen subset of $\{0,1,\dots\}$ that contains 0 missing exactly $n$
subsums is $f(n)$.

It would of course be interesting to do a similar empirical analysis for the distribution of the number of missing differences; perhaps their joint distribution could even be reduced to a simpler one using similar observations.

\section{Conjectures and open problems}\label{section.conjectures}

\begin{figure}[t]\begin{center}
\begin{picture}(380,250)
    \includegraphics{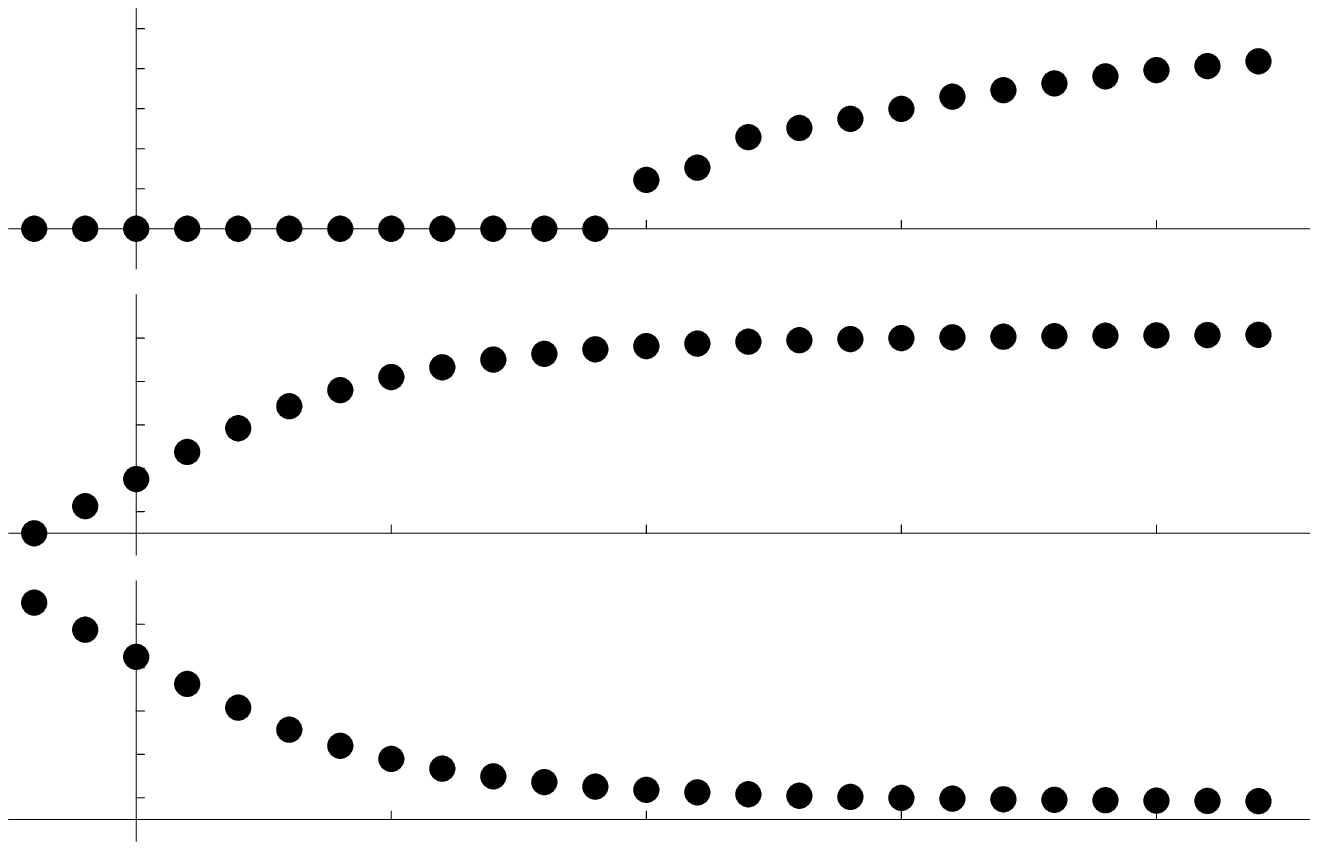}
    \put(-346,0){5}
    \put(-282,0){10}
    \put(-208,0){15}
    \put(-135,0){20}
    \put(-60,0){25}
    \put(-367,40){$0.5$}
    \put(-367,14){$0.1$}
    \put(-367,123){$0.5$}
    \put(-367,148){$0.9$}
    \put(-384,203){$0.0002$}
    \put(-384,228){$0.0004$}
    \put(-240,38){$\rho_{=}$}
    \put(-240,134){$\rho_{-}$}
    \put(-240,195){$\rho_{+}$}
    \put(-6,10){$n$}
    \put(-6,92){$n$}
    \put(-6,180){$n$}
\end{picture}\end{center}
\caption{The probability of a random subset of $\{0,\dots,n-1\}$ being \sd\ (top graph), \dd\
(middle graph), or \sdb\ (bottom graph)\label{figure.n27}}
\end{figure}

We have already conjectured, in Conjecture~\ref {main.conjecture}, that the limiting proportions of
\dd, \sdb, and \sd\ subsets of $\nset$ approach nonzero limits as $n$ tends to infinity. (As long
as the limits do in fact exist, Theorem~\ref{positive.proportion.thm} shows that they are
necessarily nonzero.) Figure~\ref{figure.n27} shows the observed proportions, for $n\le 27$, of the
subsets of $\nset$ that are \dd, \sdb, and \sd, respectively. Note particularly that each graph is
monotonic in $n$, supporting our conjecture that the limits exist. Using ten million randomly
chosen subsets of $\{0,1,\dots,999\}$, we estimate:
\begin{conjecture}
Using the notation of Conjecture~\ref {main.conjecture},
    \[
        \rho_{-}  \approx 0.93, \quad
        \rho_{+}  \approx 0.00045, \quad\text{and }
        \rho_{=} \approx 0.07.
    \]
\label{three.rho.values.conjecture}
\end{conjecture}

In fact the philosophy behind Theorem 1 suggests somewhat more: a typical subset of $\nset$ will
achieve virtually all possible sums and differences, and the ones that aren't achieved are due to
the edges of the subset. Since a positive proportion of sets have any prescribed edges, we make the
following conjecture. Define
    \begin{multline}
    \rho_{j,k} := \lim_{n\to\infty} \big( 2^{-n} \,\#\big\{ S\subset\nset :  \label {more.rho.definitions} \\
    |S+S| = 2n-1-j,\, |S-S| = 2n-1-k \big\} \big),
    \end{multline}
assuming the limit exists. Since the different set $S-S$ is symmetric about 0 and thus always has
odd cardinality, we never have $|S-S|=2n-1-k$ with $k$ odd. Therefore we conjecture:

\begin{conjecture}
For any nonnegative integers $j$ and $k$ with $k$ even, the limiting proportion $\rho_{j,k}$ defined above in~\eqref{more.rho.definitions} exists and is positive; furthermore,
\[
\sum_{j=0}^{\infty} \sum_{\substack{k=0\\k\text{ even}}}^\infty \rho_{j,k} = 1.
\]
\label{all.the.rhos.men.conjecture}
\end{conjecture}

\begin{remark}
Given Theorem \ref{theorem.average.value}, it seems reasonable to conjecture also that
\[
\sum_{j=0}^{\infty} \sum_{\substack{k=0\\k\text{ even}}}^\infty k\rho_{j,k}
    = 6 \quad\text{and}\quad \sum_{j=0}^{\infty} \sum_{\substack{k=0\\k\text{ even}}}^\infty j\rho_{j,k}
    = 10.
\]
For any particular pair $j,k$, if a single finite configuration of edges could be found that
omitted exactly $j$ possible sums and $k$ possible differences, the methods of this paper would
then show that $\rho_{j,k}>0$ (technically, that the analogous expression with
$\liminf_{n\to\infty}$ in place of $\lim_{n\to\infty}$ is positive).
\end{remark}

The last remark suggests as well the following open problem, for which a simple proof might exist, though we have not been able to find one.

\begin{conjecture}
For any nonnegative integers $j$ and $k$ with $k$ even, there exists a positive integer $n$, and a set $S\subset\nset$ with $0\in S$ and $n-1\in S$, such that $|S+S|=2n-1-j$ and $|S-S|=2n-1-k$.
\label{just.find.one.conjecture}
\end{conjecture}

\noindent Hegarty points out that his methods from \cite{hegarty} can establish both Conjecture~\ref {all.the.rhos.men.conjecture} and Conjecture~\ref {just.find.one.conjecture} in the case $j \ge k/2$.

We know~\cites{1999.Godbole.Janson.Locantore.Rapoport,2003.Nathanson} that essentially all subsets
of $\nset$ of cardinality $O(n^{1/4})$ are Sidon sets and hence \dd\ sets. More generally, we can show (perhaps in a sequel paper) that if $m=o(n^{1/2})$, then almost all subsets of $\nset$ of cardinality $m$ are \dd\ sets.

This result may indicate the presence of a threshhold. Set $p_n$ to vary with $n$, and define $n$
independent random variables $X_i$, with $X_i=1$ with probability $p_n$. This defines a random set
$A:=\{i\in\nset \colon X_i=1\}$. The observations above can then be rephrased in the following way:
if $p_n= o(n^{-1/2})$, then $A$ is a \dd\ set with probability 1 (as $n\to\infty$). We showed in
this article that if $p_n=1/2$, then $A$ is a \sd\ set with positive probability (as $n\to\infty$),
and our result is easily extended to $p_n=c>0$. An important unanswered question is ``Which
sequences $p_n$ generate a \sd\ set with positive probability?'' Perhaps our last conjecture
captures the correct notion:

\begin{conjecture}
For each $n\ge1$, let $X_{n,0},X_{n,1},\dots,X_{n,n-1}$ be independent identically distributed
random variables, and set $A_n:=\{i \colon 0\le i < n, X_{n,i}=1\}$. If both $|A_n|\to \infty$ and
$|A_n|/n \to 0$ with probability 1, then the probability that $A_n$ is \dd\ also
goes to 1.
\end{conjecture}

\end{document}